\documentclass[10pt,twocolumn,twoside]{IEEEtran}%
\usepackage{amsmath}
\usepackage{graphicx}
\usepackage{amsfonts}
\usepackage{amssymb}
\usepackage{epstopdf}%
\providecommand{\U}[1]{\protect\rule{.1in}{.1in}}

\makeatletter

\@addtoreset{theorem}{section}

\@addtoreset{corollary}{section}

\@addtoreset{lemma}{section}

\@addtoreset{definition}{section}
\makeatother
\allowdisplaybreaks
\begin{document}

\title{Optimal Routing of Energy-aware Vehicles in Networks with Inhomogeneous
Charging Nodes\thanks{The authors' work is supported in part by NSF under
Grant CNS-1139021, by AFOSR under grant FA9550-12-1-0113, by ONR under grant
N00014-09-1-1051, and by ARO under Grant W911NF-11-1-0227.}}
\author{\textbf{ S. Pourazarm, C.~G. Cassandras}\\Division of Systems Engineering and\\Center for Information and Systems Engineering, Boston University\\\texttt{sepid@bu.edu, cgc@bu.edu}}
\maketitle

\begin{abstract}
We study the routing problem for vehicles with limited energy through a
network of inhomogeneous charging nodes. This is substantially more
complicated than the homogeneous node case studied in \cite{Reno2014}. We
seek to minimize the total elapsed time for vehicles to reach their
destinations considering both traveling and recharging times at nodes when the
vehicles do not have adequate energy for the entire journey. We study two
versions of the problem. In the single vehicle routing problem, we formulate a
mixed-integer nonlinear programming (MINLP) problem and show that it can be
reduced to a lower dimensionality problem by exploiting properties of an
optimal solution. We also obtain a Linear Programming (LP) formulation
allowing us to decompose it into two simpler problems yielding near-optimal
solutions. For a multi-vehicle problem, where traffic congestion effects are
included, we use a similar approach by grouping vehicles into
\textquotedblleft subflows\textquotedblright. We also provide an alternative
flow optimization formulation leading to a computationally simpler problem
solution with minimal loss in accuracy. Numerical results are included to
illustrate these approaches.

\end{abstract}

\section{Introduction}

\label{sec1} The increasing presence of Battery-Powered Vehicles (BPVs), such
as Electric Vehicles (EVs), mobile robots and sensors, has given rise to novel
issues in classical network routing problems \cite{Laporte91}. More generally,
when entities in a network are characterized by physical attributes exhibiting
a dynamic behavior, this behavior can play an important role in the routing
decisions. In the case of BPVs, the physical attribute is \emph{energy} and
there are four BPV characteristics which are crucial in routing problems:
limited cruising range, long charge times, sparse coverage of charging
stations, and the BPV energy recuperation ability \cite{Artmeier2010} which
can be exploited. In recent years, the vehicle routing literature has been
enriched by work aiming to accommodate these BPV characteristics. For example,
by incorporating the recuperation ability of EVs, extensions to general
shortest-path algorithms are proposed in \cite{Artmeier2010} that address the
energy-optimal routing problem, with further extensions in \cite{Eisner2011}.
Charging times are incorporated into a multi-constrained optimal path planning
problem in \cite{Siddiqi2011}, which aims to minimize the length of an EV's
route and meet constraints on total traveling time, total time delay due to
signals, total recharging time and total recharging cost. In \cite{Kuller2011}%
, algorithms for several routing problems are proposed, including a
single-vehicle routing problem with inhomogeneously priced refueling stations
for which a dynamic programming based algorithm is proposed to find a least
cost path from source to destination. More recently, an EV Routing Problem
with Time Windows and recharging stations (E-VRPTW) was proposed in
\cite{TechrptMichael}, where controlling recharging times is circumvented by
simply forcing vehicles to be always fully recharged. In the Unmanned
Autonomous Vehicle (UAV) literature, a UAV routing problem with refueling
constraints is considered in \cite{Sunder2012}. A Mixed-Integer Nonlinear
Programming (MINLP) formulation is proposed with a heuristic algorithm to
determine feasible solutions.

In \cite{Reno2014} we studied the energy-constrained vehicle routing problem
in a network of homogeneous charging nodes so that the total elapsed time
(traveling time and charging time) is minimized. For a single vehicle, a MINLP
problem was formulated and, by deriving properties of the optimal solution, we
were able to decompose it into two simple Linear Programming (LP) problems.
For a multi-vehicle problem, where traffic congestion effects are included and
a system-wide objective is considered, a similar approach was used by grouping
vehicles into \textquotedblleft subflows\textquotedblright.

In this paper, we deal with the vehicle total traveling time minimization
problem in a network containing \emph{inhomogeneous} charging nodes, i.e.,
charging rates at different nodes are not identical. In fact, depending on an
outlet's voltage and current, charging an EV battery could take anywhere from
minutes to hours and the Society of Automotive Engineering (SAE) classifies
charging stations into three categories \cite{Joose2010}, \cite{Bai2010},
\cite{J1772} as shown in Tab. \ref{Tbl1}. Thus, charging rates and times are
highly dependent on the charging station class and clearly affect the solution
of our optimization problem. \begin{table}[h]
\caption{Classification of charging stations \cite{Bai2010}}%
\label{Tbl1}
\begin{center}%
\begin{tabular}
[c]{|c|c|c|c|}\hline
Charge & Nominal Supply & Max. Current & Miles per every\\
Method & Voltage(volts) & (Amps) & hour charging\\\hline
AC Level 1 & 120 VAC, 1-phase & 12 A & $<5$\\\hline
AC Level 2 & 208-240 VAC, 1-phase & 32 A & up to 62\\\hline
DC Charging & 300 - 460VDC & 400 A Max. & up to 300\\\hline
\end{tabular}
\end{center}
\end{table}As in \cite{Reno2014}, we view this as a network routing problem
where vehicles control not only their routes but also amounts to recharge at
various nodes in the network.

The contributions of this paper are as follows. For the single energy-aware
vehicle routing problem, due to the inhomogeneity in charging nodes, we can no
longer reduce the original problem to a simple LP as in \cite{Reno2014}.
However, we can still prove certain optimality properties allowing us to
reduce the dimensionality of the original problem. Further, by adopting a
locally optimal charging policy, we derive an LP formulation through which
near-optimal solutions are obtained. We note that the main difference between
this single-vehicle problem and the one considered in \cite{Kuller2011} is
that we aim to minimize the total elapsed time over a path, as opposed to a
fueling cost, and must, therefore, include two parameters per network arc:
energy consumption and traveling time. We then study a multi-vehicle
energy-aware routing problem, where a traffic flow model is used to
incorporate congestion effects. Similar to \cite{Reno2014}, by grouping
vehicles into \textquotedblleft subflows\textquotedblright\ we are able to
reduce the complexity of the original problem, although we can no longer
obtain an LP formulation. Moreover, we provide an alternative flow-based
formulation which reduces the computational complexity of the original MINLP
problem by orders of magnitude with numerical results showing little loss in optimality.

The structure of the paper is as follows. In Section \ref{sec2}, we address
the single-vehicle routing problem in a network with inhomogeneous charging
nodes and identify properties which lead to its simplification. In Section
\ref{sec3}, the multi-vehicle routing problem is formulated, first as a MINLP
and then as an alternative flow optimization problem. Simulation examples are
included illustrating our approach and providing insights on the relationship
between recharging speed and optimal routes. Conclusions and further research
directions are outlined in Section \ref{sec4}.

\section{Single Vehicle Routing}

\label{sec2} We assume that a network is defined as a directed graph
$G=(\mathcal{N},\mathcal{A})$ with $\mathcal{N}=\{1,\dots,n\}$ and
$|\mathcal{A}|=m$ . Node $i\in\mathcal{N}/\{n\}$ represents a charging station
and $(i,j)\in\mathcal{A}$ is an arc connecting node $i$ to $j$ (we assume for
simplicity that all nodes have a charging capability, although this is not
necessary). We also define $I(i)$ and $O(i)$ to be the set of start nodes
(respectively, end nodes) of arcs that are incoming to (respectively, outgoing
from) node $i$, that is, $I(i)=\{j\in\mathcal{N}|(j,i)\in\mathcal{A}\}$ and
$O(i)=\{j\in\mathcal{N}|(i,j)\in\mathcal{A}\}$.

First we deal with a single-origin-single-destination vehicle routing problem
in a network of inhomogeneous charging stations. Nodes 1 and $n$ respectively
are defined to be the origin and destination. For each arc $(i,j)\in
\mathcal{A}$, there are two cost parameters: the required traveling time
$\tau_{ij}$ and the energy consumption $e_{ij}$. Note that $\tau_{ij}>0$ (if
nodes $i$ and $j$ are not connected, then $\tau_{ij}=\infty$), whereas
$e_{ij}$ is allowed to be negative due to a BPV's potential energy
recuperation effect \cite{Artmeier2010}. Letting the vehicle's charge capacity
be $B$, we assume that $e_{ij}<B$ for all $(i,j)\in\mathcal{A}$. Since we are
considering a single vehicle's behavior, we assume that it will not affect the
overall network's traffic state, therefore, $\tau_{ij}$ and $e_{ij}$ are fixed
depending on given traffic conditions at the time the single-vehicle routing
problem is solved. Clearly, this cannot apply to the multi-vehicle case in the
next section, where the decisions of multiple vehicle routes affect traffic
conditions, thus influencing traveling times and energy consumption. Since the
BPV has limited battery energy it may not be able to reach the destination
without recharging. Thus, recharging amounts at charging nodes $i\in
\mathcal{N}$ are also decision variables.

We denote the selection of arc $(i,j)$ and energy recharging amount at node
$i$ by $x_{ij}\in\{0,1\}$, $i,j\in\mathcal{N}$ and $r_{i}\geq0$,
$i\in\mathcal{N}/\{n\}$, respectively. We also define $g_{i}$ as the charging
time per unit of energy for charging node $i$, i.e., the reciprocal of a
charging rate for each node. Without loss of generality we assume $g_{n}=0$.
Moreover, we use $E_{i}$ to represent the vehicle's residual battery energy at
node $i$. Then, for all $E_{j},\,j\in O(i)$, we have:
\begin{equation}
E_{j}=\left\{
\begin{array}
[c]{ll}%
E_{i}+r_{i}-e_{ij} & \text{if }x_{ij}=1\\
0 & \text{otherwise}%
\end{array}
\right.  \label{ExpandEj}%
\end{equation}
which can also be expressed as
\[
E_{j}=\sum_{i\in I(j)}(E_{i}+r_{i}-e_{ij})x_{ij},\quad x_{ij}\in\{0,1\}
\]
The single vehicle's objective is to determine a path from $1$ to $n$, as well
as recharging amounts, so as to minimize the total elapsed time to reach the
destination. We formulate this as a Mixed Integer Nonlinear Programming
(MINLP) problem:
\begin{gather}
\min_{x_{ij},r_{i},\,\,i,j\in\mathcal{N}}\quad\sum_{i=1}^{n}\sum_{j=1}^{n}%
\tau_{ij}x_{ij}+\sum_{i=1}^{n}\sum_{j=1}^{n}r_{i}g_{i}x_{ij}\label{obj}\\
s.t.\quad\sum_{j\in O(i)}x_{ij}-\sum_{j\in I(i)}x_{ji}=b_{i},\quad\text{for
each }i\in\mathcal{N}\label{flowConv}\\
b_{1}=1,\,b_{n}=-1,\,b_{i}=0,\text{ for }i\neq1,n\label{bi}\\
E_{j}=\sum_{i\in I(j)}(E_{i}+r_{i}-e_{ij})x_{ij},\text{ for }j=2,\dots
,n\label{EiConv}\\
0\leq E_{i}\leq B,\quad E_{1}\text{ given},\text{ for each }i\in
\mathcal{N}\label{Ei}\\
x_{ij}\in\{0,1\},\quad r_{i}\geq0 \label{controls}%
\end{gather}
We will refer to this problem as \textbf{P1}. The constraints (\ref{flowConv}%
)-(\ref{bi}) stand for the flow conservation, which implies that only one path
starting from node $i$ can be selected, i.e., $\sum_{j\in O(i)}x_{ij}\leq1$.
It is easy to check that this also implies $x_{ij}\leq1$ for all $i,j$ since
$b_{1}=1$, $I(1)=\varnothing$. Constraint (\ref{EiConv}) represents the
vehicle's energy dynamics where the only nonlinearity in this formulation
appears. Finally, (\ref{Ei}) indicates that the vehicle cannot run out of
energy before reaching a node or exceed a given capacity $B$. All other
parameters are predetermined according to the network topology. A crucial
difference between \textbf{P1} and the MINLP introduced in \cite{Reno2014} is
that here the charging rates $g_{i}$ in (\ref{obj}) are node-dependent.

\subsection{Properties}

Rather than directly tackling the MINLP problem, we derive some key properties
of an optimal solution which will enable us to reduce \textbf{P1} to a
lower-dimension problem. In particular, there are $m+2(n-1)$ decision
variables in \textbf{P1} (because of the $E_{j}$ variables in (\ref{EiConv})),
which we will show how to reduce to $m+(n-1)$. The main difficulty in this
problem lies in the coupling of the decision variables $x_{ij}$ and $r_{i}$ in
(\ref{EiConv}) and the following lemma will enable us to eliminate $r_{i}$
from (\ref{obj}).\newline\textbf{Lemma 1: }Given (\ref{obj})-(\ref{controls}),
an optimal solution $\{x_{ij},r_{i}\},\,\,i,j\in\mathcal{N}$ satisfies:
\textbf{ }
\begin{gather}
\sum_{i=1}^{n}\sum_{j=1}^{n}(r_{i}x_{ij}-e_{ij}x_{ij})g_{i}=\sum_{i=1}^{n}%
\sum_{j=1}^{n}(E_{j}-E_{i})g_{i}x_{ij}\label{Lemma1}\\
=\sum_{i=1}^{n}\sum_{j=1}^{n}E_{j}(g_{i}-g_{j})x_{ij}-E_{1}g_{1}
\label{Lemma1b}%
\end{gather}
\emph{Proof}: Multiplying both sides of (\ref{ExpandEj}) by $g_{i}$ gives:
\[
E_{j}g_{i}=%
\begin{cases}
(E_{i}+r_{i}-e_{ij})g_{i} & \text{if }x_{ij}=1,\\
0 & \text{otherwise }.
\end{cases}
\]
which can be expressed as $\sum_{i\in I(j)}E_{j}g_{i}x_{ij}=\sum_{i\in
I(j)}(E_{i}+r_{i}-e_{ij})g_{i}x_{ij}$. Summing both sides over $j=2,\ldots,n$
and rearranging yields:
\begin{align*}
&  \sum_{j=2}^{n}\sum_{i\in I(j)}E_{j}g_{i}x_{ij}-\sum_{j=2}^{n}\sum_{i\in
I(j)}E_{i}g_{i}x_{ij}\\
&  =\sum_{j=2}^{n}\sum_{i\in I(j)}(r_{i}-e_{ij})g_{i}x_{ij}%
\end{align*}
Based on (\ref{ExpandEj}), $E_{i}=0$ for all nodes which are not in the
selected path. Thus we can rewrite the equation above as%
\[
\sum_{i=1}^{n}\sum_{j=1}^{n}(r_{i}x_{ij}-e_{ij}x_{ij})g_{i}=\sum_{i=1}^{n}%
\sum_{j=1}^{n}(E_{j}-E_{i})g_{i}x_{ij}%
\]
which establishes (\ref{Lemma1}). Finally, (\ref{Lemma1b}) follows by
observing that if $P$ is an optimal path we can re-index nodes so that
$P=\{1,...,n\}$ with $g_{n}=0$. Thus, we have $\sum_{i=1}^{n}\sum_{j=1}%
^{n}E_{i}g_{i}x_{ij}=E_{1}g_{1}+\ldots+E_{n-1}g_{n-1}$ which can also be
written as $E_{1}g_{1}+\sum_{i=2}^{n}\sum_{j=2}^{n}E_{j}g_{j}x_{ij}$ where
$x_{ij}=0$ for all $(i,j)$ not in the optimal path. Therefore,%
\[
\sum_{i=1}^{n}\sum_{j=1}^{n}(E_{j}-E_{i})g_{i}x_{ij}=\sum_{i=1}^{n}\sum
_{j=1}^{n}E_{j}(g_{i}-g_{j})x_{ij}-E_{1}g_{1}%
\]
which proves (\ref{Lemma1b}).$\blacksquare$\newline\textbf{Lemma 2:} If
$\sum_{i}r_{i}^{\ast}>0$ in the optimal routing policy, then $E_{n}^{\ast}=0$.
\newline\emph{Proof:} This is the same as the homogeneous charging node case;
see Lemma 2 in \cite{Reno2014}.\newline In view of Lemma 1, we can replace
$\sum_{i=1}^{n}\sum_{j=1}^{n}r_{i}g_{i}x_{ij}$ in (\ref{obj}) and eliminate
the presence of $r_{i}$, $i=2,\ldots,n-1$, from the objective function and the
constraints. This results in the following MINLP problem referred to as
\textbf{P2}:
\begin{gather}
\displaystyle{\min_{\substack{x_{ij},E_{i}\\i,j\in\mathcal{N}}}\sum_{i=1}%
^{n}\sum_{j=1}^{n}\big(\tau_{ij}x_{ij}+e_{ij}g_{i}x_{ij}+E_{j}(g_{i}%
-g_{j})x_{ij}\big)-E_{1}g_{1}}\label{PIIobj}\\
s.t\quad\sum_{j\in O(i)}x_{ij}-\sum_{j\in I(i)}x_{ji}=b_{i}\\
b_{1}=1,\ b_{n}=-1,\ b_{i}=0\quad\text{for}\quad i\neq1,n\\
0\leqslant E_{j}-(E_{i}-e_{ij})x_{ij}\leqslant B\quad\forall\ i,j\in
N\label{PIIConst}\\
0\leqslant E_{i}\leqslant B,\ E_{1}\text{ given},\ \forall\ i\in N\\
x_{ij}\in\{0,1\}
\end{gather}
This new formulation has only $m+(n-1)$ decision variables compared to
$m+2(n-1)$ in \textbf{P1}. Constraint (\ref{PIIConst}) is derived from
(\ref{EiConv}). Assuming $x_{ij}=1$, i.e. arc $(i,j)$ is part of the optimal
path, we can recover $r_{i}=E_{j}-E_{i}+e_{ij}$ and Constraint (\ref{PIIConst}%
) is added to prevent any vehicle from exceeding its capacity $B$ in an
optimal path. Solving this problem gives both an optimal path and residual
battery energy at each node.

Although \textbf{P2} has fewer decision variables, it is still a MINLP which
is hard to solve for large networks. Specifically, the CPU time is highly
dependent on the number of nodes and arcs in the network. In what follows we
introduce a \emph{locally optimal} charging policy, leading to a simpler
problem, by arguing as follows. Looking at (\ref{PIIobj}), the term
${\sum_{i=1}^{n}\sum_{j=1}^{n}E_{j}(g_{i}-g_{j})x_{ij}}$ is minimized by
selecting each $E_{j}$ depending on the sign of ${(g_{i}-g_{j})}$%
:\newline\textbf{Case 1}: $g_{i}-g_{j}<0$, i.e., node $i$ has a faster
charging rate than node $j$. Therefore, $E_{j}$ should get its maximum
possible value, which is $B-e_{ij}$. This implies that the vehicle must be
maximally charged at node $i$.\newline\textbf{Case 2}: $g_{i}-g_{j}\geqslant
0$, i.e., node $j$ has a faster or same charging rate as node $i$. In this
case, $E_{j}$ should get its minimum value $E_{j}=0$. This implies that the
vehicle should get the minimum charge needed at node $i$ in order to reach node $j$.\newline We
define $\mathbf{\pi_{C}}$ to be the charging policy specified as above and
note that it does not guarantee the global optimality of $E_{i}$ thus selected
in (\ref{PIIobj}) which can easily be checked by a counterexample. However, it allows us to decompose
the optimal routing problem from the optimal charging problem. If, in
addition, we consider only solutions for which the vehicle is recharged at
least once (otherwise, the vehicle is not energy-constrained and the problem
is of limited interest), we can obtain the following result.\newline%
\textbf{Theorem 1: } If $\sum_{i}r_{i}^{\ast}>0$ (i.e. the vehicle has to be
recharged at least once), then under charging policy $\mathbf{\pi_{C}}$, the
solution $x_{ij}^{\ast}$, $i,j\in\mathcal{N}$, of the original problem
(\ref{obj}) can be determined by solving the LP problem:%
\begin{gather}
\min_{x_{ij,}i,j\in\mathcal{N}}\sum_{i=1}^{n}\sum_{j=1}^{n}\big(\tau
_{ij}+e_{ij}g_{i}+K(g_{i}-g_{j})\big)x_{ij}\label{PIIIobj}\\
K=%
\begin{cases}
B-e_{ij} & \text{if }g_{i}<g_{j},\\
0 & \text{otherwise }.
\end{cases}
\label{PIIIK}\\
s.t.\quad\sum_{j\in O(i)}x_{ij}-\sum_{j\in I(i)}x_{ji}=b_{i}\\
b_{1}=1,b_{n}=-1,b_{i}=0\hspace{12pt}fori\neq1,n\\
0\leqslant x_{ij}\leqslant1
\end{gather}
\textit{Proof}: Applying charging policy $\mathbf{\pi_{C}}$ in (\ref{PIIobj})
we can change objective function to ${\sum_{i=1}^{n}\sum_{j=1}^{n}%
\big(\tau_{ij}+e_{ij}g_{i}+K(g_{i}-g_{j})\big)x_{ij}-E_{1}g_{1}}$ where $K$ is
as in (\ref{PIIIK}). Therefore, $x_{ij}^{\ast}$ can be determined by the
following problem:
\begin{gather}
\min_{x_{ij,}i,j\in\mathcal{N}}\sum_{i=1}^{n}\sum_{j=1}^{n}\big(\tau
_{ij}+e_{ij}g_{i}+K(g_{i}-g_{j})\big)x_{ij}-E_{1}g_{1}\nonumber\\
K=%
\begin{cases}
B-e_{ij} & \text{if }g_{i}<g_{j},\\
0 & \text{otherwise }.
\end{cases}
\nonumber\\
s.t.\quad\sum_{j\in O(i)}x_{ij}-\sum_{j\in I(i)}x_{ji}=b_{i},\quad\text{for
each }i\in\mathcal{N}\nonumber\\
b_{1}=1,\,b_{n}=-1,\,b_{i}=0,\text{ for }i\neq1,n\nonumber\\
x_{ij}\in\{0,1\}\nonumber
\end{gather}
which is a typical shortest path problem formulation. Moreover, according to
the property of minimum cost flow problems \cite{Hillier}, the above integer
programming problem is equivalent to the LP with the integer restriction on
$x_{ij}$ relaxed. Finally, since $E_{1}$ and $g_{1}$ are given, the problem
reduces to (\ref{PIIIobj}), which proves the theorem. $\blacksquare$%
\newline\textbf{Remark 1}. If $g_{i}=g_{j}$ for all $i,j$ in (\ref{PIIIobj}),
the problem reduces to the homogeneous charging node case studied in
\cite{Reno2014} with the same optimal LP formulation as in Theorem 1. With
$g_{i}\neq g_{j}$ however, the LP formulation cannot guarantee global
optimality, although the routes obtained through Theorem 1 may indeed be
optimal (see Section II.C), in which case the optimal charging amounts are
obtained as described next.\newline
\vspace{-5mm}
\subsection{Determination of optimal recharging amounts $r_{i}^{\ast}$}

Once we determine an optimal route $P$, it is relatively easy to find a
feasible solution for $r_{i}$, $i\in P$, to satisfy the constraint
(\ref{EiConv}) and minimize the total charging time on the selected path. It
is obvious that the optimal charging amounts $r_{i}^{\ast}$ are non-unique in
general. Without loss of generality we re-index nodes so that we may write
$P=\{1,...,n\}$. Then, the problem resulting in an optimal charging policy is
\begin{gather}
\min_{r_{i},\text{ }i\in P}\quad\sum_{i\in P}g_{i}r_{i}\label{obj4}\\
s.t.\quad E_{i+1}=E_{i}+r_{i}-e_{i,i+1}\nonumber\\
0\leq E_{i}\leq B,\quad E_{1}\text{ given}\nonumber\\
r_{i}\geq0\text{ for all }i\in\mathcal{N}\nonumber
\end{gather}
This is an LP where $E_{i}$ and $r_{i}$ are decision variables. Unlike the
homogeneous charging node problem in \cite{Reno2014} where the objective
function includes charging prices $p_{i}$ associated with nodes, i.e.,
$\sum_{i\in P}p_{i}r_{i}$, this is not the case here, since there is a
tradeoff between selecting faster-charging nodes and possible higher costs at
such nodes. However, the advantage of the decoupling approach is that if an
optimal path is determined, an additional cost minimization problem can be
formulated to determine optimal charging times at nodes on this path.

\begin{figure}[h]
\begin{center}
\includegraphics[scale=0.35]{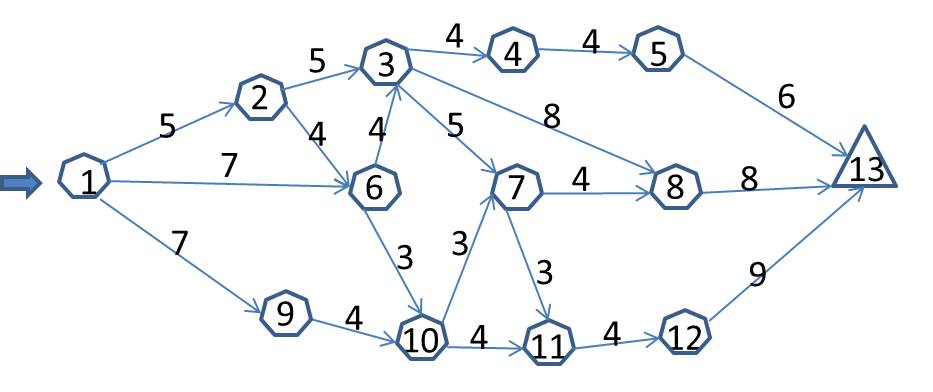}
\end{center}
\caption{13-node network with inhomogeneous charging nodes}\label{NetSVR}
\end{figure}
\vspace{-5mm}
\subsection{Numerical Example}

We consider a 13-node network as shown in Fig.\ref{NetSVR} where the travel
distance of each arc is as shown. For simplicity here we assume $\tau
_{ij}=e_{ij}=d_{ij}$ and solve the problem for different configurations of
charging stations in the network. The optimal paths obtained by solving
different formulations are shown in Tab.\ref{TblExSVR1} where $G=[g_{1}%
,\ g_{2},\ ,...,g_{n-1}]$. We can see that all formulations result in the
\emph{same} optimal path while their computational complexities are
drastically different (from around 250 sec for \textbf{P1} to less than 2 sec
for the LP formulation). Once the optimal path is determined, we can easily
solve (\ref{obj4}) to determine optimal charging amounts as well.
\begin{table}[th]
\caption{Optimal paths obtained by solving Problems \textbf{P-1}, \textbf{P-2}
and LP for different charging node configurations}%
\label{TblExSVR1}
\begin{center}%
\begin{tabular}
[c]{|c|c|c|}\hline
G & Problem & Path\\\hline
\lbrack1 0.2 1 0.1 1 0.2 0.1 1 1 1 1 1] & \textbf{P1}\& \textbf{P2}: &
1$\rightarrow$2$\rightarrow$3$\rightarrow$4$\rightarrow$5$\rightarrow$13\\
& \textbf{LP}: & 1$\rightarrow$2$\rightarrow$3$\rightarrow$4$\rightarrow
$5$\rightarrow$13\\\hline
\lbrack1 1 1 1 1 1 0.2 0.1 1 0.1 1 1] & \textbf{P1}\& \textbf{P2}: &
1$\rightarrow$6$\rightarrow$10$\rightarrow$7$\rightarrow$8$\rightarrow$13\\
& \textbf{LP}: & 1$\rightarrow$6$\rightarrow$10$\rightarrow$7$\rightarrow
$8$\rightarrow$13\\\hline
\lbrack1 1 1 1 1 1 1 1 1 0.1 1 1] & \textbf{P1}\& \textbf{P2}: &
1$\rightarrow$6$\rightarrow$10$\rightarrow$7$\rightarrow$8$\rightarrow$13\\
& \textbf{LP}: & 1$\rightarrow$6$\rightarrow$10$\rightarrow$7$\rightarrow
$8$\rightarrow$13\\\hline
\lbrack1 1 1 1 1 1 1 1 1 1 1 1] & \textbf{P1}\& \textbf{P2}: & 1$\rightarrow
$2$\rightarrow$3$\rightarrow$4$\rightarrow$5$\rightarrow$13\\
& \textbf{LP}: & 1$\rightarrow$2$\rightarrow$3$\rightarrow$4$\rightarrow
$5$\rightarrow$13\\\hline
\lbrack1 1 1 1 1 1 1 0.1 0.1 1 1 1] & \textbf{P1}\& \textbf{P2}: &
1$\rightarrow$9$\rightarrow$10$\rightarrow$7$\rightarrow$8$\rightarrow$13\\
& \textbf{LP}: & 1$\rightarrow$9$\rightarrow$10$\rightarrow$7$\rightarrow
$8$\rightarrow$13\\\hline
\end{tabular}
\end{center}
\end{table}

\section{Multiple Vehicle Routing}

\label{sec3} We now investigate the multi-vehicle routing problem in a network
with inhomogeneous charging nodes, where we seek to optimize a
\emph{system-wide objective} by routing and charging vehicles. The main
technical difficulty in this case is that we need to consider the influence of
traffic congestion on both traveling time and energy consumption. A second
difficulty is that of implementing an optimal routing policy. In the case of a
centrally controlled system consisting of mobile robots, sensors or any type
of autonomous vehicles this can be accomplished through appropriately
communicated commands. In the case of vehicles with individual drivers,
implementation requires signaling mechanisms and possibly incentive structures
to enforce desired routes assigned to vehicles, bringing up a number of
additional research issues. In the sequel, we limit ourselves to resolving the
first difficulty before addressing implementation challenges.

If we proceed as in the single vehicle case, i.e., determining a path
selection through $x_{ij}^{k}$, $i,j\in\mathcal{N}$, and recharging amounts
$r_{i}^{k}$, $i\in\mathcal{N}/\{n\}$ for all vehicles $k=1,\ldots,K$, for some
$K$, then the dimensionality of the solution space is prohibitive. Moreover,
the inclusion of traffic congestion effects introduces additional
nonlinearities in the dependence of the travel time $\tau_{ij}$ and energy
consumption $e_{ij}$ on the traffic flow through arc $(i,j)$, which now depend
on $x_{ij}^{1},\cdots,x_{ij}^{K}$. Instead, as in \cite{Reno2014}, we proceed
by grouping subsets of vehicles into $N$ \textquotedblleft
subflows\textquotedblright\ where $N$ may be selected to render the problem
manageable (we will discuss the effect of $N$ in Section III.C).

Let all vehicles enter the network at node 1 and let $R$ denote the rate of
vehicles arriving at this node. Viewing vehicles as defining a \emph{flow}, we
divide them into $N$ \emph{subflows} each of which may be selected so as to
include the same type of homogeneous vehicles (e.g., large vehicles vs smaller
ones or vehicles with the same initial energy). Thus, all vehicles in the same
subflow follow the same routing and recharging decisions so that we only
consider control at the subflow level rather than individual vehicles. Note
that asymptotically, as $N\rightarrow\infty$, we can recover routing at the
individual vehicle level. Clearly, not all vehicles in our system are BPVs,
therefore, not part of our optimization process. These can be treated as
uncontrollable interfering traffic for our purposes and can be readily
accommodated in our analysis, as long as their flow rates are known. However,
for simplicity, we will assume here that every arriving vehicle is a BPV and
joins a subflow. Our objective is to determine optimal routes and energy
recharging amounts for each vehicle subflow so as to minimize the total
elapsed time of these flows from origin to destination. The decision variables
are $x_{ij}^{k}\in\{0,1\}$ and $r_{i}^{k}$ for all arcs $(i,j)$ and subflows
$k=1,\ldots,N$. Given traffic congestion effects, the time and energy
consumption on each arc depends on the values of $x_{ij}^{k}$ and the fraction
of the total flow rate $R$ associated with each subflow $k$; the simplest such
flow allocation is one where each subflow is assigned $R/N$. Let
$\mathbf{x_{ij}}=(x_{ij}^{1},\cdots,x_{ij}^{N})^{T}$ and $\mathbf{r_{i}%
}=(r_{i}^{1},\cdots,r_{i}^{N})^{T}$. Then, we denote the traveling time and
corresponding energy consumption of the $k$th vehicle subflow on arc $(i,j)$
by $\tau_{ij}^{k}(\mathbf{x_{ij}})$ and $e_{ij}^{k}(\mathbf{x_{ij}})$
respectively. As already mentioned, $\tau_{ij}^{k}(\mathbf{x_{ij}})$ and
$e_{ij}^{k}(\mathbf{x_{ij}})$ can also incorporate the influence of
uncontrollable (non-BPV) vehicle flows, which can be treated as parameters in
these functions. Similar to the single vehicle case, we use $E_{i}^{k}$ to
represent the residual energy of subflow $k$ at node $i$, given by the
aggregated residual energy of all vehicles in the subflow. If the subflow does
not go through node $i$, then $E_{i}^{k}=0$. Similar to \cite{Reno2014}, the
problem is formulated as follows:
\begin{gather}
\min_{\mathbf{x_{ij}},\mathbf{r_{i}},\,\,i,j\in\mathcal{N}}\quad\sum_{i=1}%
^{n}\sum_{j=1}^{n}\sum_{k=1}^{N}\left(  \tau_{ij}^{k}(\mathbf{x_{ij}}%
)+r_{i}^{k}g_{i}x_{ij}^{k}\right) \label{objM}\\
\text{s.t.}\quad\sum_{j\in O(i)}x_{ij}^{k}-\sum_{j\in I(i)}x_{ji}^{k}%
=b_{i},\quad\text{for each }i\in\mathcal{N}\label{flowConvM}\\
b_{1}=1,\,b_{n}=-1,\,b_{i}=0,\text{ for }i\neq1,n\label{biM}\\
E_{j}^{k}=\sum_{i\in I(j)}(E_{i}^{k}+r_{i}^{k}-e_{ij}^{k}(\mathbf{x_{ij}%
}))x_{ij}^{k},\quad j=2,\dots,n\label{EiConvM}\\
E_{1}^{k}\text{ is given},\quad E_{i}^{k}\geq0,\quad\text{ for each }%
i\in\mathcal{N}\label{EiM}\\
x_{ij}^{k}\in\{0,1\},\quad r_{i}^{k}\geq0 \label{controlsM}%
\end{gather}
We will refer to this problem \textbf{P3}. The difference from the MINLP
formulated in \cite{Reno2014} is that we consider different charging rates
$g_{i}$ in the objective function. \textbf{P3 }contains $N(m+2(n-1))$ decision
variables and is difficult to solve. However, as in the single-vehicle case,
we are able to establish some properties allowing us to simplify it.
\vspace{-5mm}
\subsection{Properties}

Even though the term $\tau_{ij}^{k}(\mathbf{x_{ij}})$ in the objective
function is no longer linear in general, for each subflow $k$ the constraints
(\ref{flowConvM})-(\ref{controlsM}) are still similar to the single-vehicle
case. Consequently, we can derive similar useful properties in the form of the
following lemmas (proofs are very similar to those of the single-vehicle case
and are omitted).\newline\textbf{Lemma 3: } An optimal solution
$\{\mathbf{x_{ij}},\mathbf{r_{i}}\},\,\,i,j\in\mathcal{N}$satisfies:
\begin{align}
&  \sum_{i=1}^{n}\sum_{j=1}^{n}\sum_{k=1}^{N}r_{i}^{k}g_{i}x_{ij}^{k}%
-\sum_{i=1}^{n}\sum_{j=1}^{n}\sum_{k=1}^{N}e_{ij}^{k}(\mathbf{x_{ij}}%
)g_{i}x_{ij}^{k}\label{Lemma2a}\\
&  =\sum_{i=1}^{n}\sum_{j=1}^{n}\sum_{k=1}^{N}(E_{j}^{k}-E_{i}^{k})g_{i}%
x_{ij}^{k}\nonumber\\
&  =\sum_{i=1}^{n}\sum_{j=1}^{n}\sum_{k=1}^{N}E_{j}^{k}(g_{i}-g_{j})x_{ij}%
^{k}-\sum_{k=1}^{N}E_{1}^{k}g_{1} \label{Lemma2b}%
\end{align}
\textbf{Lemma 4:} If $\sum_{i}r_{i}^{k\ast}>0$ in the optimal routing policy,
then $E_{n}^{k\ast}=0$ for $k=1,...,N$. \newline In view of Lemma 3, we can
replace ${\sum_{i=1}^{n}\sum_{j=1}^{n}\sum_{k=1}^{N}r_{i}^{k}g_{i}x_{ij}^{k}}$
in (\ref{objM}) through (\ref{Lemma2b}) and eliminate $r_{i}^{k}$,
$i=1,\ldots,n-1$, $k=1,\ldots,N$, from the objective function (\ref{objM}).
The term $\sum_{k=1}^{N}E_{1}^{k}g_{1}$ is also removed because it has a fixed
value. Thus, we introduce a new MINLP formulation to determine $x_{ij}^{k\ast
}$ and $E_{i}^{k\ast}$ for all $i,j\in\mathcal{N}$ and $k=1,\ldots,N$ as
follows:
\begin{gather}
\min_{\substack{x_{ij}^{k},E_{i}^{k}\\i,j\in\mathcal{N}}}\sum_{i=1}^{n}%
\sum_{j=1}^{n}\sum_{k=1}^{N}\big(\tau_{ij}^{k}(\mathbf{x_{ij}})+(e_{ij}%
^{k}(\mathbf{x_{ij}})g_{i}+E_{j}^{k}(g_{i}-g_{j}))x_{ij}^{k}%
\big)\label{PIVobj}\\
\text{s.t.}\quad\sum_{j\in O(i)}x_{ij}^{k}-\sum_{j\in I(i)}x_{ji}^{k}=b_{i}\\
b_{1}=1,\ b_{n}=-1,\ b_{i}=0\ \quad\text{for}\ i\neq1,n\nonumber\\
0\leq E_{j}^{k}-(E_{i}^{k}-e_{ij}^{k}(\mathbf{x_{ij}}))x_{ij}^{k}\leq
B^{k}\quad\mbox{for each }(i,j)\in A\label{PIVconst}\\
E_{1}^{k}\text{is given},\quad E_{i}^{k}\geq0,\quad\text{for each }i\in
N\label{PIVconst2}\\
x_{ij}^{k}\in\{0,1\}
\end{gather}
We call this problem \textbf{P4}. Note that inequality (\ref{PIVconst}) is
derived from (\ref{EiConvM}). Assuming $x_{ij}^{k}=1$, i.e., arc $(i,j)$ is
part of the optimal path for the $k$th subflow, $r_{i}^{k}=e_{ij}%
^{k}(\mathbf{x_{ij}})+E_{j}^{k}-E_{i}^{k}$. Thus, (\ref{PIVconst}) ensures the
optimal solution $E_{i}^{k\ast}$ results in a feasible charging amount for the
$k$th subflow, $0\leq r_{i}^{k}\leq B^{k}$ where $B^{k}$ is the maximum
charging amount $k$th subflow can get. This value should be predetermined for
each subflow based on the vehicle types and the fraction of total inflow in
it. Similar to \textbf{P2} in the single-vehicle case, once we determine
$E_{i}^{k\ast}$ we can simply calculate optimal charging amounts using
(\ref{EiConvM}). Although \textbf{P4} has $N(m+(n-1))$ decision variables,
which is fewer than \textbf{P3}, its complexity still highly depends on the
network size and number of subflows. Similar to the charging policy
$\mathbf{\pi_{C}}$ used in Theorem 1, we introduce a charging policy by
arguing as follows. Looking at (\ref{PIVobj}), the term ${\sum_{i=1}^{n}%
\sum_{j=1}^{n}\sum_{k=1}^{N}E_{j}^{k}(g_{i}-g_{j})x_{ij}^{k}}$ is minimized by
selecting each $E_{j}^{k}$ depending on the sign of ${(g_{i}-g_{j})}$:
\newline\textbf{Case 1}: $g_{i}<g_{j}$, i.e., node $i$ has faster charging
rate than node $j$. Therefore, $E_{j}^{k}$ should get its maximum value, i.e.,
the $k$th subflow should get its maximum charge at node $i$.\newline%
\textbf{Case 2}: $g_{i}\geqslant g_{j}$, i.e., the charging rate of node $j$
is greater than or equal to node $i$. Therefore, $E_{j}^{k}$ should get its
minimum value of $0$. This implies that the $k$th subflow should get the
minimum charge needed at node $i$ in order to reach node $j$.\newline Applying this policy in
(\ref{PIVobj}) and changing the objective function accordingly we introduce
problem \textbf{P5} as follows:
\begin{gather}
\min_{x_{ij}^{k}}\sum_{i=1}^{n}\sum_{j=1}^{n}\sum_{k=1}^{N}\big[\tau_{ij}%
^{k}(\mathbf{x_{ij}})+(e_{ij}^{k}(\mathbf{x_{ij}})g_{i}+K(g_{i}-g_{j}%
))x_{ij}^{k}\big]\label{PVobj}\\
K=%
\begin{cases}
B^{k}-e_{ij}^{k}(\mathbf{x_{ij}}) & \text{if }g_{i}<g_{j},\label{PVK}\\
0 & \text{{}}\mbox{otherwise}
\end{cases}
\\
\text{ s.t.}\quad\sum_{j\in O(i)}x_{ij}^{k}-\sum_{j\in I(i)}x_{ji}^{k}=b_{i}\\
b_{1}=1,b_{n}=-1,b_{i}=0\quad\text{for }i\neq1,n\nonumber\\
x_{ij}^{k}\in\{0,1\}
\end{gather}
Unlike the single-vehicle case, the objective function is no longer
necessarily linear in $x_{ij}^{k}$, therefore, (\ref{PVobj}) cannot be further
simplified into an LP problem as in Theorem 1. The computational effort
required to solve this problem with $Nm$ decision variables, depends on the
dimensionality of the network and the number of subflows. Nonetheless, from
the transformed formulation above, we are still able to separate the
determination of routing variables $x_{ij}^{k}$ from recharging amounts
$r_{i}^{k}$. Similar to the single-vehicle case, once the routes are
determined, we can obtain $r_{i}^{k}$ satisfying the energy constraints
(\ref{EiConvM})-(\ref{EiM}) while minimizing $\sum_{k=1}^{N}\sum_{i\in P^k}%
r_{i}^{k}g_{i}$ where $P^K$ is the optimal path of the $k$th subflow. Next, we present an alternative formulation of
(\ref{objM})-(\ref{controlsM}) leading to a computationally simpler solution
approach.\newline\textbf{Remark 2}. If $g_{i}=g_{j}$ for all $i,j$ in
(\ref{PVobj}), the problem reduces to the homogeneous charging node case with
the exact same MINLP formulation as in \cite{Reno2014} for obtaining an
optimal path. However, \textbf{P5} cannot guarantee an optimal solution
because of the locally optimal charging policy $\mathbf{\pi_{C}}$ which may
not be feasible in a globally optimal solution $(x_{ij}^{k\ast},E_{i}^{k\ast
})$.

\subsection{Flow control formulation}

We begin by relaxing the binary variables in (\ref{controlsM}) by letting
$0\leq x_{ij}^{k}\leq1$. Thus, we switch our attention from determining a
single path for any subflow $k$ to several possible paths by treating
$x_{ij}^{k}$ as the normalized vehicle flow on arc $(i,j)$ for the $k$th
subflow. This is in line with many network routing algorithms in which
fractions $x_{ij}$ of entities are routed from a node $i$ to a neighboring
node $j$ using appropriate schemes ensuring that, in the long term, the
fraction of entities routed on $(i,j)$ is indeed $x_{ij}$. Following this
relaxation, the objective function in (\ref{objM}) is changed to:
\[
\min_{\mathbf{x_{ij}},\mathbf{r_{i}},\,\,i,j\in\mathcal{N}}\quad\sum_{i=1}%
^{n}\sum_{j=1}^{n}\sum_{k=1}^{N}\tau_{ij}^{k}(\mathbf{x_{ij}})+\sum_{i=1}%
^{n}\sum_{k=1}^{N}r_{i}^{k}g_{i}%
\]
Moreover, the energy constraint (\ref{EiConvM}) needs to be adjusted
accordingly. Let $E_{ij}^{k}$ represent the fraction of residual energy of
subflow $k$ associated with the $x_{ij}^{k}$ portion of the vehicle flow
exiting node $i$. Therefore, constraint (\ref{EiM}) becomes $E_{ij}^{k}\geq0$.
We can now capture the relationship between the energy associated with subflow
$k$ and the vehicle flow as follows:
\begin{equation}
\left[  \sum_{h\in I(i)}(E_{hi}^{k}-e_{hi}^{k}(\mathbf{x_{ij}}))+r_{i}%
^{k}\right]  \cdot\frac{x_{ij}^{k}}{\sum_{h\in I(i)}x_{hi}^{k}}=E_{ij}%
^{k}\label{EiConvM2}%
\end{equation}%
\begin{equation}
\frac{E_{ij}^{k}}{\sum_{j\in O(i)}E_{ij}^{k}}=\frac{x_{ij}^{k}}{\sum_{j\in
O(i)}x_{ij}^{k}}\label{same}%
\end{equation}
In (\ref{EiConvM2}), the energy values of different vehicle flows entering
node $i$ are aggregated and the energy corresponding to each portion exiting a
node, $E_{ij}^{k}$ , $j\in O(i)$, is proportional to the corresponding
fraction of vehicle flows, as expressed in (\ref{same}). Clearly, this
aggregation of energy leads to an approximation, since one specific vehicle
flow may need to be recharged in order to reach the next node in its path,
whereas another might have enough energy without being recharged. This
approximation foregoes controlling recharging amounts at the individual
vehicle level and leads to approximate solutions of the original problem
(\ref{objM})-(\ref{controlsM}). Several numerically based comparisons are
provided in the next section showing little or no loss of optimality relative
to the solution of (\ref{objM}). Adopting this formulation with
$x_{ij}^{k}\in\lbrack0,1]$ instead of $x_{ij}^{k}\in\{0,1\}$, we obtain the
following simpler nonlinear programming problem (NLP):
\begin{gather}
\min_{\mathbf{x_{ij}},\mathbf{r_{i}},\,\,i,j\in\mathcal{N}}\quad\sum_{i=1}%
^{n}\sum_{j=1}^{n}\sum_{k=1}^{N}\tau_{ij}^{k}(\mathbf{x_{ij}})+\sum_{i=1}%
^{n}\sum_{k=1}^{N}r_{i}^{k}g_{i}\label{objM3}\\
\text{s.t.}\quad\sum_{j\in O(i)}x_{ij}^{k}-\sum_{j\in I(i)}x_{ji}^{k}%
=b_{i},\quad\text{for each }i\in\mathcal{N}\label{flowConvM3}\\
b_{1}=1,\,b_{n}=-1,\,b_{i}=0,\text{ for }i\neq1,n\nonumber\\
\left[  \sum_{h\in I(i)}(E_{hi}^{k}-e_{hi}^{k}(\mathbf{x_{ij}}))+r_{i}%
^{k}\right]  \cdot\frac{x_{ij}^{k}}{\sum_{h\in I(i)}x_{hi}^{k}}=E_{ij}%
^{k}\label{EiConvM3}\\
\frac{E_{ij}^{k}}{\sum_{j\in O(i)}E_{ij}^{k}}=\frac{x_{ij}^{k}}{\sum_{j\in
O(i)}x_{ij}^{k}}\label{sameM3}\\
E_{ij}^{k}\geq0,\label{EiM3}\\
0\leq x_{ij}^{k}\leq1,\quad r_{i}^{k}\geq0\label{controlsM3}%
\end{gather}
As in our previous analysis, we are able to eliminate $\mathbf{r_{i}}$ from
the objective function in (\ref{objM3}) as follows.\newline\textbf{Lemma 5: }
For each subflow $k=1,\ldots,N$,
\begin{align}
&  \sum_{i=1}^{n}r_{i}^{k}g_{i}=\nonumber\\
&  \sum_{i=1}^{n}\sum_{j=1}^{n}e_{ij}^{k}(\mathbf{x_{ij}})g_{i}+\sum_{i=1}%
^{n}\sum_{j\in O(i)}E_{ij}^{k}g_{i}-\sum_{i=1}^{n}\sum_{h\in I(i)}E_{hi}%
^{k}g_{i}\nonumber\\
&  =\sum_{i=1}^{n}\sum_{j=1}^{n}e_{ij}^{k}(\mathbf{x_{ij}})g_{i}+\sum
_{i=1}^{n}\sum_{j\in O(i)}E_{ij}^{k}(g_{i}-g_{j})\nonumber
\end{align}
\emph{Proof}: Multiplying (\ref{EiConvM3}) by $g_{i}$ and summing over all
$i=1,\ldots,n$ , then using (\ref{flowConvM3}) and (\ref{sameM3}) proves the
lemma. $\blacksquare$\newline Using Lemma 5 we change the objective function
(\ref{objM3}) to:
\begin{equation}
\sum_{i=1}^{n}\sum_{j=1}^{n}\sum_{k=1}^{N}\big(\tau_{ij}^{k}(\mathbf{x_{ij}%
})+e_{ij}^{k}(\mathbf{x_{ij}})g_{i}\big)+\sum_{i=1}^{n}\sum_{j=1}^{n}%
\sum_{k=1}^{N}E_{ij}^{k}(g_{i}-g_{j})\label{objtoNLP}\\
\end{equation}
Once again, we adopt a charging policy $\mathbf{\pi_{C}}$ as follows:\newline%
\textbf{Case 1}: If $g_{i}<g_{j}$, then $E_{ij}^{k}$ gets its maximum value
$(B^{k}-e_{ij}^{k}(\mathbf{x_{ij}}))x_{ij}^{k}$.\newline\textbf{Case 2}: If
$g_{i}\geq g_{j}$, then $E_{ij}^{k}$ gets its minimum value $0$.\newline
Applying this policy in (\ref{objtoNLP}) we can transform the objective
function (\ref{objM3}) to (\ref{objM4}) and determine near-optimal routes
$x_{ij}^{k\ast}$ by solving the following NLP:
\begin{gather}
\min_{\substack{\mathbf{x_{ij}}\\i,j\in\mathcal{N}}}\quad\sum_{k=1}^{N}%
\sum_{i=1}^{n}\sum_{j=1}^{n}\left[  \tau_{ij}^{k}(\mathbf{x_{ij}})+e_{ij}%
^{k}(\mathbf{x_{ij}})g_{i}+K(g_{i}-g_{j})\right]  \label{objM4}\\
K=%
\begin{cases}
(B^{k}-e_{ij}^{k}(\mathbf{x_{ij}}))x_{ij}^{k} & \text{if }g_{i}<g_{j}%
,\nonumber\\
0 & \text{{}otherwise}%
\end{cases}
\\
\text{s.t.}\quad\sum_{j\in O(i)}x_{ij}^{k}-\sum_{j\in I(i)}x_{ji}^{k}%
=b_{i},\quad\text{for each }i\in\mathcal{N}\nonumber\\
b_{1}=1,\,b_{n}=-1,\,b_{i}=0,\text{ for }i\neq1,n\nonumber\\
0\leq x_{ij}^{k}\leq1\nonumber
\end{gather}
Once again, there is no guarantee of global optimality. The values of $r_{i}^{k}$,
$i=1,\ldots,n$, $k=1,\ldots,N$, can then be determined so as to satisfy the
energy constraints (\ref{EiConvM3})-(\ref{EiM3}), and minimizing $\sum
_{k=1}^{N}\sum_{i \in P^k}^{n}r_{i}^{k}g_{i}$. Finally, if $g_{i}=g_{j}$ for all
$i,j$ in (\ref{objM4}), the problem reduces to the homogeneous charging node
case with the same exact NLP flow control formulation as in \cite{Reno2014}.
\vspace{-5mm}
\subsection{Numerical Examples}

\begin{figure}[ptbh]
\begin{center}
\includegraphics[scale=0.4]{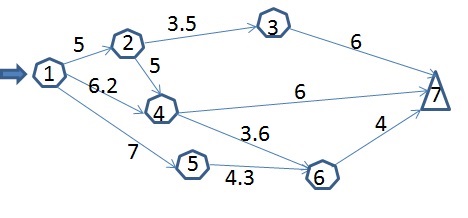}
\end{center}
\caption{A 7-node network with inhomogeneous charging nodes.}%
\label{smallGr}%
\end{figure}

We consider a specific example which includes traffic congestion and energy
consumption functions. The relationship between the speed and density of a
vehicle flow is typically estimated as follows (see \cite{Ioannou96}):
\begin{equation}
v(k(t))=v_{f}\bigg(1-\left(  \frac{k(t)}{k_{jam}}\right)  ^{p}\bigg)^{q}%
\label{VvsD}%
\end{equation}
where $v_{f}$ is the reference speed on the road without traffic, $k(t)$
represents the density of vehicles on the road at time $t$ and $k_{jam}$ the
saturated density for a traffic jam. The parameters $p$ and $q$ are
empirically identified for actual traffic flows. In our multi-vehicle routing
problem, we are interested in the relationship between the density of the
vehicle flow and traveling time on an arc $(i,j)$, i.e., $\tau_{ij}%
^{k}(\mathbf{x_{ij}})$. Given a network topology (i.e., a road map), the
distances $d_{ij}$ between nodes are known. Moreover, we do not include
uncontrollable vehicle flows in our example for simplicity. In our approach,
we need to identify $N$ subflows and we do so by evenly dividing the entire
vehicle inflow into $N$ subflows, each of which has $R/N$ vehicles per unit
time. Thus, $k_{jam}$ in this case can be set as $N$, implying that we do not
want all vehicles to go through the same path, hence the the arc $(i,j)$
density is $\sum_{k}x_{ij}^{k}$. Therefore, the time subflow $k$ spends on arc
$(i,j)$ becomes
\[
\tau_{ij}^{k}(\mathbf{x_{ij}})=\big(d_{ij}\cdot x_{ij}^{k}\cdot\frac{R}%
{N}\big)\big(v_{f}(1-(\frac{\sum_{k}x_{ij}^{k}}{N})^{p})^{q}\big)^{-1}%
\]
As for $e_{ij}^{k}(\mathbf{x_{ij}})$, we assume the energy consumption rates
of subflows on arc $(i,j)$ are all identical, proportional to the distance
between nodes $i$ and $j$, giving
\[
e_{ij}^{k}(\mathbf{x_{ij}})=e\cdot d_{ij}\cdot\frac{R}{N}%
\]
In order to verify the accuracy of different formulations, we numerically solve the optimal and near-optimal problems \textbf{P4} and  \textbf{P5}. For the latter, (\ref{PVobj}) 
becomes:
\begin{gather}
\min_{\substack{x_{ij}^{k}\\i,j\in\mathcal{N}}}\quad\sum_{i=1}^{n}\sum
_{j=1}^{n}\sum_{k=1}^{N}\big[\frac{d_{ij}x_{ij}^{k}\frac{R}{N}}{v_{f}%
(1-(\frac{\sum_{k}x_{ij}^{k}}{N})^{p})^{q}}+eg_{i}d_{ij}\frac{R}{N}x_{ij}%
^{k}+\nonumber\\
K(g_{i}-g_{j})\big]\label{objSp}\\
K=%
\begin{cases}
(B^{k}-ed_{ij}\frac{R}{N})x_{ij}^{k} & \text{if }g_{i}<g_{j},\label{PVK}\\
0 & \text{{}otherwise}%
\end{cases}
\\
\text{ s.t. for each }k\in\{1,...,N\}:\nonumber\\
\sum_{j\in O(i)}x_{ij}^{k}-\sum_{j\in I(i)}x_{ji}^{k}=b_{i}\\
b_{1}=1,b_{n}=-1,b_{i}=0\quad\text{for }i\neq1,n\nonumber\\
x_{ij}^{k}\in\{0,1\}
\end{gather}
For simplicity, we let $v_{f}=1$ mile/min, $R=1$ vehicle/min, $p=2,\,\,q=2$
and $e=1$. The network topology used is that of Fig.\ref{smallGr}, with the
distance of each arc as shown. \begin{table}[ptbh]
\caption{Numerical results for sample problem}%
\label{table2}
\begin{center}%
\begin{tabular}
[c]{|c|c|c|}\hline
& \textbf{P4} & \textbf{P5}\\\hline
N & 2 & 2\\\hline
obj & 34.459444 & 35.0444\\\hline
routes & $1\rightarrow5\rightarrow6\rightarrow7$ & $1\rightarrow
5\rightarrow6\rightarrow7$\\
& $1\rightarrow4\rightarrow7$ & $1\rightarrow4\rightarrow7$\\\hline\hline
N & 6 & 6\\\hline
obj & 29.228750 & 29.6187\\\hline
routes & $1\rightarrow2\rightarrow3\rightarrow7(\times2)$ & $1\rightarrow
2\rightarrow3\rightarrow7(\times2)$\\
& $1\rightarrow4\rightarrow7(\times2)$ & $1\rightarrow4\rightarrow7(\times
2)$\\
& $1\rightarrow5\rightarrow6\rightarrow7(\times2)$ & $1\rightarrow
5\rightarrow6\rightarrow7(\times2)$\\\hline\hline
N & 10 & 10\\\hline
obj & 29.290898 & 29.6419\\\hline
routes & $1\rightarrow2\rightarrow3\rightarrow7(\times3)$ & $1\rightarrow
2\rightarrow3\rightarrow7(\times3)$\\
& $1\rightarrow4\rightarrow7(\times4)$ & $1\rightarrow4\rightarrow7(\times
4)$\\
& $1\rightarrow5\rightarrow6\rightarrow7(\times3)$ & $1\rightarrow
5\rightarrow6\rightarrow7(\times3)$\\\hline\hline
N & 18 & 18\\\hline
obj &  & 29.542209\\\hline
routes &  & $1\rightarrow2\rightarrow3\rightarrow7(\times5)$\\
& Solver Error & $1\rightarrow4\rightarrow7(\times7)$\\
&  & $1\rightarrow5\rightarrow6\rightarrow7(\times6)$\\\hline\hline
N & 30 & 30\\\hline
obj &  & 29.501931\\\hline
routes &  & $1\rightarrow2\rightarrow3\rightarrow7(\times9)$\\
& Solver Error & $1\rightarrow4\rightarrow7(\times11)$\\
&  & $1\rightarrow5\rightarrow6\rightarrow7(\times10)$\\\hline
\end{tabular}
\end{center}
\end{table}Tab. \begin{figure}[ptbh]
\begin{center}
\includegraphics[scale=0.52]{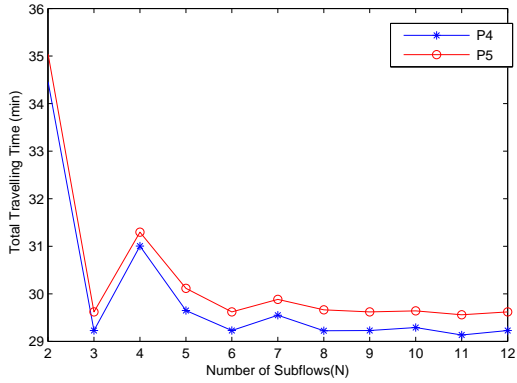}
\end{center}
\caption{Performance as a function of $N$ (No. of subflows)}%
\label{NvsObj}%
\end{figure}\ref{table2} shows both optimal routes and locally optimal routes
obtained by solving \textbf{P4} and \textbf{P5} respectively for different
values of $N\in\lbrack1,\ldots,30]$ and $G=[1\ 1\ 1\ 1\ 0.1\ 1]$. We observe
that vehicles are mainly distributed through three routes and the traffic
congestion effect makes the flow distribution differ from the shortest path.
The number of decision variables (hence, the solution search space) rapidly
increases with the number of subflows. However, looking at Fig. \ref{NvsObj}
which gives the performance in terms of our objective functions in
(\ref{PIVobj}) and (\ref{objSp}) as a function of the number of subflows,
observe that the optimal objective value (\textbf{P4}) quickly converges
around $N=8$. Thus, even though the best solution is found when $N=11$, a
near-optimal solution can be determined under a small number of subflows. This
suggests that one can rapidly approximate the asymptotic solution of the
multi-vehicle problem (dealing with individual vehicles routed so as to
optimize a system-wide objective) based on a relatively small value of $N$.
Another observation is that although \textbf{P5} is a suboptimal formulation
it results in the same paths as those obtained by solving \textbf{P4}. Moreover,
the cost difference is small over different $N$.

Next, we obtain a solution to the same problem (\ref{objSp}) using the NLP
formulation (\ref{objM4}) with $0\leq x_{ij}^{k}\leq1$. Since in this example
all subflows are identical, we can further combine all $x_{ij}^{k}$ over each
$(i,j)$, leading to the $N$-subflow relaxed problem:
\begin{gather}
\min_{\substack{x_{ij}^{k}\\i,j\in\mathcal{N}}}\sum_{i=1}^{n}\sum_{j=1}%
^{n}\left[  \frac{d_{ij}x_{ij}R}{v_{f}(1-(x_{ij})^{p})^{q}}+eg_{i}%
d_{ij}Rx_{ij}+K(g_{i}-g_{j})\right]  \label{objSp2}\\
K=%
\begin{cases}
(B-ed_{ij}R)x_{ij} & \text{if }g_{i}<g_{j},\label{PVK}\\
0 & \text{{}otherwise}%
\end{cases}
\nonumber\\
s.t.\quad\sum_{j\in O(i)}x_{ij}-\sum_{j\in I(i)}x_{ji}=b_{i},\quad\text{for
each }i\in\mathcal{N}\nonumber\\
b_{1}=1,\,b_{n}=-1,\,b_{i}=0,\text{ for }i\neq1,n\nonumber\\
0\leq x_{ij}\leq1\nonumber
\end{gather}
This is a relatively easy to solve NLP problem. Using the same parameter
settings as before, we obtain the objective value of $28.5645$ mins and the
optimal routes are: $35.938\%$ of vehicle flow: $(1\rightarrow4\rightarrow7)$;
$28.605\%$ of vehicle flow: $(1\rightarrow2\rightarrow3\rightarrow7)$;
$35.457\%$ of vehicle flow: $(1\rightarrow5\rightarrow6\rightarrow7)$.
Compared to the best solution ($N=11$) in Fig. \ref{NvsObj}, the difference in
objective values between the integer and flow-based solutions is less than
$2\%$. This supports the effectiveness of a solution based on a limited number
of subflows in the MINLP problem.

\textbf{CPU time Comparison}. Tab. \ref{table3} compares the computational
effort in terms of CPU time for problems \textbf{P3}, \textbf{P5} and the flow
control formulation to find optimal routes for the sample network shown in
Fig. \ref{smallGr}. Our results show that the flow control formulation results
in a reduction of about 4 orders of magnitude in CPU time with approximately the same
objective function value.\begin{table}[ptbh]
\caption{CPU time for sample problem}%
\label{table3}
\begin{center}%
\begin{tabular}
[c]{c|c|c|c}\hline
\textbf{Fig.\ref{smallGr} Net.} & \textbf{P3} & \textbf{P5} & NLP
approx.\\\hline
N & 3(near opt) & 3(near opt) & -\\\hline
obj & 29.2287 & 29.6187 & 28.5645\\\hline
CPU time(sec) & 17122.644 & 190.7711 & 1.4\\\hline
\end{tabular}
\end{center}
\end{table}


\section{Conclusions and future work}

\label{sec4} We have studied the problem of minimizing the total elapsed time
for energy-constrained vehicles to reach their destinations, including
recharging when there is no adequate energy for the entire journey. In
contrast to our earlier work \cite{Reno2014}, we have considered inhomogeneous
charging rates at nodes. For a single vehicle, we have shown how to reduce the
complexity of this problem. For a multi-vehicle problem, where traffic
congestion effects are considered, we used a similar approach by aggregating
vehicles into subflows and seeking optimal routing decisions for each such
subflow. We also developed an alternative flow-based formulation which yields
approximate solutions with a computational cost reduction of several orders of
magnitude, so they can be used in problems of large dimensionality. Numerical
examples show these solutions to be near-optimal. We have also found that a
low number of subflows is adequate to obtain convergence to near-optimal solutions.

\bibliographystyle{IEEEtran}
\bibliography{IEEEabrv,renowang}

\end{document}